\documentclass[a4paper,11pt]{amsart}
\usepackage{eucal}
\usepackage{cite}
\usepackage{amsmath,amsthm,amssymb}
\usepackage{amsfonts}
\usepackage{latexsym}
\usepackage{iftex}
\ifpdf\else
  \ifXeTeX\else
    \PassOptionsToPackage{dvipdfmx}{hyperref}
  \fi
\fi
\usepackage{hyperref}
\usepackage{eucal}
\usepackage{mathrsfs}
\usepackage{txfonts}
\theoremstyle{plain}
\newtheorem{theorem}{Theorem}[section]
\newtheorem*{thm*}{Theorem}

\newtheorem{proposition}[theorem]{Proposition}

\newtheorem{lemma}[theorem]{Lemma}
\newtheorem{corollary}[theorem]{Corollary}

\theoremstyle{definition}
\newtheorem{definition}[theorem]{Definition}

\makeatletter
\renewenvironment{proof}[1][\proofname]{\par
  \normalfont
  \topsep6\p@\@plus6\p@ \trivlist
  \item[\hskip\labelsep{\bfseries #1}\@addpunct{\bfseries.}]\ignorespaces
}{%
  \endtrivlist
}
\renewcommand{\proofname}{proof}
\theoremstyle{remark}
\newtheorem{remark}[theorem]{Remark}

\numberwithin{equation}{section}

\title{Nilpotent approximation and completion of $\E_\infty$-algebra objects of stable symmetric monoidal model categories}
\date{\today} 
\author{Yuki Kato}
\thanks{The author was supported by Grants-in-Aid for Scientific Research No.~23K03080, Japan Society for the Promotion of Science.}
\address{National institute of technology, Kurume college, 
	      1-1-1, Komorino, Kurume, Fukuoka, JAPAN 830-8555.}
\email{kato{\_}051@kurume-nct.ac.jp}

\newcommand{\A}{\mathbb{A}}
\newcommand{\Z}{\mathbb{Z}}
\newcommand{\K}{\mathbb{K}}
\newcommand{\E}{\mathbb{E}}

\newcommand{\proj}{\mathbb{P}}
\newcommand{\Spec}{\mathrm{Spec}}

\newcommand{\Hom}{\mathrm{Hom}}

\newcommand{\Map}{\mathrm{Map}}

\newcommand{\Mod}{\mathrm{Mod}}
\newcommand{\CAlg}{\mathrm{CAlg}}

\newcommand{\Fun}{\mathrm{Fun}}

\newcommand{\SShv}{\mathbf{S}.\,\mathbf{Shv}}
\newcommand{\Sm}{\mathbf{Sm}}

\newcommand{\tm}{\tilde{\mathfrak{m}}}
\newcommand{\FMod}{\mathrm{Mod}^{\rm fi}}

\newcommand{\AMod}{\mathrm{AMod}}
\newcommand{\APMod}{\mathrm{APMod}}

\usepackage{enumerate}
\usepackage{array}
\usepackage[all]{xy} 
\usepackage{delarray}
\def\qed{{\hfill $\Box$}}
\begin{document}
\thispagestyle{empty}
\begin{abstract}
Motivated by derived completion in homotopy theory~\cite{Goodwillie1,Goodwillie2,Goodwillie3,HA,zbMATH06622452}, this article develops a nilpotent approximation theory for Smith ideals, extending adic completion for commutative rings to monoid objects in locally presentable symmetric monoidal abelian categories and to stable symmetric monoidal model categories. Our principal result is the formal completeness theorem (Theorem~\ref{TheoremA1}): finite generation of a Smith ideal forces completeness of its nilpotent approximation. This gives a precise categorical analogue of~\cite{compFG}, while remaining distinct from classical adic completion.

As applications, we obtain an almost-mathematics version of completion theory and prove that homotopical nilpotent approximation is complete for weakly compact Smith ideals.

As a geometric application of the general theory, we study the nilpotent approximation determined by the canonical morphism $\mathbf{MGL}\to\K$ in Voevodsky's symmetric monoidal category of motivic spectra. This yields a natural $K$-theoretic approximation of algebraic cobordism; we prove that it is homotopically complete and satisfies Bott periodicity, and we establish a mod-$l$ Gabber rigidity theorem for the corresponding approximation of $\mathbf{MGL}/l$ by $\K/l$.
\end{abstract}

\maketitle
\section{Introduction}
\label{sec-Introduction} 
Let $A$ be a commutative ring, and let $I$ be an ideal of $A$. Then the filtration $\{I^n\}_{n \ge 1}$ defines a topological ring structure on $A$, called the {\it $I$-adic topology}. The projective limit $\hat{A}=\varprojlim A/I^n$ is called the {\it completion} of $A$ with respect to the $I$-adic topology. For each $m \ge 1$, write $\widehat{I^m}=\varprojlim_n I^m/I^{m+n}$. In general, $\hat{A}$ is complete for the topology defined by the induced filtration $\{ \widehat{I^n} \}_{n \ge 1}$, but not for the $\hat{I}$-adic topology $\{ \hat{I}^{\,n} \}_{ n\ge 1}$. An example is given in Yekutieli~\cite[Example 1.8]{compFG}. For this reason, we use the terminology ``nilpotent approximation'' rather than ``completion'' in this paper. For any finitely generated ideal $I$, Yekutieli~\cite[Corollary 3.6]{compFG} shows that all induced homomorphisms $I^n \hat{A} \to \hat{I}^n \to \widehat{I^n}$ are isomorphisms for each $n \ge 1$, implying that $\hat{A}$ is complete for the $\hat{I}$-adic topology.

Our approach is motivated by completion in homotopy theory, especially by Goodwillie calculus and Preygel's completion of $\E_\infty$-rings~\cite{Goodwillie1,Goodwillie2,Goodwillie3,HA,zbMATH06622452}. This paper's nilpotent approximation of Smith ideals and classical adic completion of quotient rings are conceptually distinct constructions.  We prove a categorical analogue of ~\cite[Corollary 3.6]{compFG}: in the Smith-ideal tower, finite generation implies completeness by a formal argument in the arrow-category framework.

To fix notation, we briefly recall Goodwillie calculus and the Taylor approximation theory of functors, following Lurie~\cite{HA}. Let $\mathcal{C}$ be an $\infty$-category. For any $n \ge 0$, an $n$-cube is a functor from the nerve $N(\mathbf{P}([n]))$ of the power set $\mathbf{P}([n])$ of the totally ordered set $[n]=\{0,\,1,\, \ldots,\,n\}$ to $\mathcal{C}$, where morphisms in $\mathbf{P}([n])$ are given by inclusions of subsets of $[n]$.

\begin{definition}[\cite{HA} p.1016, Definition 6.1.1.6]
Let $\mathcal{C}$ be an $\infty$-category. We say that $\mathcal{C}$ is
{\it differentiable} if $\mathcal{C}$ satisfies the following
conditions:
\begin{enumerate}
\item The $\infty$-category $\mathcal{C}$ admits finite limits.
\item The $\infty$-category $\mathcal{C}$ admits all sequential colimits. That is, for any diagram $X :N( \Z_{\geq 0} )\to
      \mathcal{C}$, the colimit of $X$ exists in $\mathcal{C}$.
\item The colimit functor $\varinjlim: \Fun( N( \Z_{\geq 0}),
       \,\mathcal{C} ) \to \mathcal{C} $ is left exact.
\end{enumerate}
\end{definition}

\begin{definition}[\cite{HA} p.1015, Definition 6.1.1.3]
Let $F:\mathcal{C} \to \mathcal{D}$ be a functor between
$\infty$-categories.  Suppose that $\mathcal{C}$ admits finite colimits
and $\mathcal{D}$ admits finite limits.  For any $n \ge 0$, we say that $F :\mathcal{C} \to  \mathcal{D} $ is {\it $n$-excisive} if $F$ carries strongly coCartesian $n$-cubes in $\mathcal{C}$ to Cartesian $n$-cubes in $\mathcal{D}$ (see Lurie~\cite[p.1015, Definition 6.1.1.2]{HA}), and we let
$\mathrm{Exc}^n(\mathcal{C},\,\mathcal{D})$ denote the full subcategory
of $\Fun( \mathcal{C},\,\mathcal{D})$ spanned by $n$-excisive functors.
\end{definition}

\begin{theorem}[\cite{HA} p.1016, Theorem 6.1.1.10]
\label{Taylor-tower} Let $\mathcal{C}$ and $\mathcal{D}$ be
$\infty$-categories. Suppose that $\mathcal{C}$ admits finite colimits
and has a final object, and $\mathcal{D}$ is differentiable.
Then the inclusion functor
\[
\mathrm{Exc}^n(\mathcal{C},\,\mathcal{D}) \to \Fun(\mathcal{C},\,\mathcal{D})
\]
admits a left adjoint
\[
P_n: \Fun(\mathcal{C},\,\mathcal{D}) \to \mathrm{Exc}^n(\mathcal{C},\,\mathcal{D}).
\]
Furthermore, the left adjoint $P_n$ is left exact. \qed
\end{theorem}
Let $F: \mathcal{C} \to \mathcal{D}$ be a functor between
$\infty$-categories. Suppose that $\mathcal{C}$ and $\mathcal{D}$
satisfy the hypotheses of Theorem~\ref{Taylor-tower}. Then we obtain a
tower of functors
\[
\cdots \to P_n(F) \to P_{n-1}(F) \to \cdots \to P_1(F) \to P_0(F)
\]
from $\mathcal{C}$ to $\mathcal{D}$. We call the tower $\{P_n(F)\}_{n \ge 0}$ the {\it Taylor tower} of $F$. Goodwillie calculus is the theory of approximating functors by $n$-excisive functors.

Preygel~\cite{zbMATH06622452} defines the completion of a morphism $f: R \to R'$ of $\E_\infty$-rings by the homotopy limit
\[
   \hat{R} =  \varprojlim P^n ( \mathrm{id}_{\CAlg_{/R'}} ) (f) \in  \CAlg_{/R'}.
\]
However, explicitly constructing an $n$-excisive approximation from the definition of Goodwillie calculus is often complicated. Inspired by Heuts's work~\cite{Heuts2018} on Goodwillie calculus for $\infty$-categories and by Hovey's Smith ideal theory~\cite{Smith-ideals}, we introduce a more elementary approximation theory for $\E_\infty$-algebra objects in stable symmetric monoidal model categories, using nilpotent Smith ideals of degree $n+1$ in place of $n$-excisive approximations. This generalizes adic completion in commutative algebra with respect to ideals. As an application, we establish an almost-mathematics version of adic completion, which we call the ``almost completion'' theory, without assuming that the defining ideal is tight~\cite[p. 32, Definition 5.1.5]{GR}.

Our principal contribution is the formal completeness theorem (Theorem~\ref{TheoremA1}): in the Smith-ideal nilpotent tower, finite generation of a Smith ideal forces completeness of its nilpotent approximation. This yields a precise analogue of Yekutieli's finite-generation completeness phenomenon~\cite{compFG}, while remaining distinct from classical adic completion.

We then generalize nilpotent approximation to homotopy theory and prove that, for weakly compact Smith ideals, it gives a homotopical analogue of completeness.

As a geometric application of this general framework, rather than asking directly which $K$-theoretic properties algebraic cobordism itself should satisfy, we study the nilpotent approximation determined by the canonical morphism $\mathbf{MGL}\to\K$. This yields a natural approximation of algebraic cobordism by algebraic $K$-theory that is homotopically complete and satisfies Bott periodicity. In addition, we prove a mod-$l$ Gabber rigidity statement for the corresponding approximation of $\mathbf{MGL}/l$ by $\K/l$~\cite{Gabber-rigid}.

This paper is organized as follows. Section~\ref{sec:finite} recalls finitely generated objects in categories. In Section~\ref{sec:Smith}, we define the nilpotent approximation theory of Smith ideals, which generalizes completion theory for commutative algebras. The main result of Section~\ref{sec:Smith} is the completeness theorem for finitely generated Smith ideals (Theorem~\ref{TheoremA1}). As an application, the final part of Section~\ref{sec:Smith} introduces almost nilpotent approximation (completion) for almost modules; unlike Gabber--Ramero's completion theory, our construction does not rely on tight ideals. We also prove an almost analogue of Gabber--Ramero~\cite[p.147, Theorem 5.3.24]{GR}. In Section~\ref{sec:model}, we generalize nilpotent approximation to homotopy theory: we define homotopically nilpotent approximation for Smith ideals in stable symmetric monoidal model categories and prove that, for weakly compact Smith ideals, it gives a homotopical analogue of completeness. Finally, Section~\ref{sec:MSP} gives a geometric application of the general theory: it constructs the approximation of algebraic cobordism by algebraic $K$-theory, proves its homotopical completeness and Bott periodicity, and then derives a mod-$l$ Gabber rigidity theorem for the corresponding approximation of $\mathbf{MGL}/l$ by $\K/l$.

Throughout the paper, we fix an uncountable regular cardinal $\kappa$. Thus ``small'' sets and categories mean $\kappa$-small sets and categories. We work with locally $\kappa$-presentable categories, proper combinatorial model categories, and closed symmetric monoidal structures, that is, symmetric monoidal structures equipped with internal Hom objects.


\section{Finitely generated (presented) objects of the arrow categories of locally presentable symmetric monoidal abelian categories}
\label{sec:finite} 
We recall the notions of finitely generated and finitely presented objects in a category:
\begin{definition}
Let $\mathcal{C}$ be a locally presentable category. An object $C$ of $\mathcal{C}$ is said to be {\it finitely generated} (resp. {\it finitely presented}) if, for every filtered inductive system $\{ M_{\alpha} \}$, the induced map
\[
  \varinjlim \Hom_\mathcal{C}(C,\,  M_{\alpha} ) \to   \Hom_\mathcal{C}(C,\, \varinjlim \ M_{\alpha} )
\]
is injective (resp. bijective).
\end{definition}

\begin{proposition}[{\rm c.f.} \cite{LocPresentable} p.53, Proposition 1.69]
\label{f-g}
Let $\mathcal{C}$ be a locally presentable pointed category, and let $C$ be an object of $\mathcal{C}$. Then the following conditions are equivalent:
\begin{enumerate}[(1)]
 \item The object $C$ is finitely generated.
 \item For any filtered inductive system $\{M_\alpha\}$ in $\mathcal{C}$ whose transition morphisms are monomorphisms, the induced map
\[
  \varinjlim \Hom_\mathcal{C}(C,\, M_\alpha) \to \Hom_\mathcal{C}(C,\, \varinjlim M_\alpha)
\]
is bijective.
\item There is an epimorphism $p:F \to C$ from a finitely presented object $F$.
\end{enumerate}
\end{proposition}
\begin{proof}
Assume that $C$ is finitely generated. We prove condition~(2). Given a filtered inductive system $\{M_\alpha\}$ with monic transition morphisms, let $Q_\alpha$ denote the cokernel of the canonical morphism $M_\alpha \to \varinjlim M_\alpha$ for each $\alpha$. Then the filtered colimit of $\{Q_\alpha\}$ is zero. Since $C$ is finitely generated, the induced map
\[
  \varinjlim \Hom_\mathcal{C}(C,\, Q_\alpha) \to \Hom_\mathcal{C}(C,\, \varinjlim Q_\alpha)=0
\]
is injective. Hence the colimit $\varinjlim \Hom_\mathcal{C}(C,\, Q_\alpha)$ is also zero, which proves condition~(2).

Next, assume that condition~(2) holds. Since $\mathcal{C}$ is locally presentable, $C$ can be expressed as a filtered colimit of a filtered system $\{F_\alpha\}$ of finitely presented objects. Any locally presentable category is complete, well-powered, and admits (strong-epi, mono)-factorizations. For each $\alpha$, let $C_\alpha$ denote the image of the morphism $F_\alpha \to C$. Condition~(2) implies that the identity $\mathrm{id}_C:C\to C$ factors through some $C_\alpha$. Hence, for such $\alpha$, the monomorphism $C_\alpha \to C$ is an isomorphism.

Finally, assume condition~(3). Since the induced map $\Hom_\mathcal{C}(C,\,M) \to \Hom_\mathcal{C}(F,\,M)$ is injective for every object $M$, it is immediate that $C$ is finitely generated. \qed
\end{proof}

In this paper, we are mainly interested in symmetric monoidal categories generated by a projective monoidal unit. We therefore recall the notion of a generator:
\begin{definition}[\cite{zbMATH01216133} p.127, Chapter V, Section 7]
Let $\mathcal{C}$ be a category. A collection $\mathcal{S}$ of objects of $\mathcal{C}$ is said to be a {\it generator} of $\mathcal{C}$ if, whenever $f \neq g: M \to N$ are parallel morphisms, there exist an object $S \in \mathcal{S}$ and a morphism $s:S \to M$ such that $f \circ s \neq g \circ s$.
\end{definition}

\begin{lemma}
\label{generator}
Let $\mathcal{C}$ be an additive category generated by an object $V$, and let $0$ denote the zero object. Then $M=0$ if and only if the set $\Hom_{\mathcal{C}}(V,\,M)$ consists of a single point.
\end{lemma}
\begin{proof}
The only-if direction is clear. Assume that $M \neq 0$ and $\Hom_\mathcal{C}(V,\,M)=\{*\}$. Then there exists an object $N$ such that $\Hom_\mathcal{C}(N,\,M)\neq \{*\}$. Let $f:N \to M$ be a nonzero morphism. Since $V$ generates $\mathcal{C}$, there exists a morphism $s:V \to N$ such that $f \circ s:V \to M$ is nonzero. This contradicts the assumption that $\Hom_\mathcal{C}(V,\,M)=\{*\}$. \qed
\end{proof}

Let $\mathcal{C}$ be an abelian category. An object $P$ of $\mathcal{C}$ is said to be {\it projective} if the functor $\Hom_\mathcal{C}(P,-)$ is exact.

\begin{proposition}
\label{closed-compact} Let $(\mathcal{C},\,\otimes )$ be a locally presentable
closed symmetric monoidal abelian category. Assume that
$\mathcal{C}$ is generated by a finitely presented projective monoidal unit $V$. Then the classes of finitely generated and finitely presented objects are closed under tensor products.
\end{proposition}
\begin{proof}

Let $\Map(-,\,-):\mathcal{C}^{\rm op} \times \mathcal{C} \to \mathcal{C}$ denote the internal Hom functor right adjoint to the monoidal product $-\otimes-$. Let $M$ and $N$ be finitely presented objects, let $\{X_\alpha\}$ be a $\kappa$-filtered system of objects of $\mathcal{C}$, and let $X$ denote its filtered colimit.

The canonical morphisms $\{\psi_\alpha:X_\alpha \to X\}$ induce a morphism
\[
 \Psi: \varinjlim  \Map(N,\,X_\alpha) \to \Map(N,\,X).
\]
By adjunction, the unitor $\eta_N:V \otimes N \to N$ induces a natural bijection
\[
  \Hom_\mathcal{C}(V,\,\Map(N,\,Y)) \simeq \Hom_\mathcal{C}(N,\,Y)
\]
for every object $Y$ of $\mathcal{C}$. Since $V$ and $N$ are finitely presented, it follows that
\[
\Psi_*(V):\Hom_\mathcal{C} \left(V,\,\varinjlim  \Map(N,\,X_\alpha)\right) \to \Hom_\mathcal{C} \left(V,\,\Map(N,\, X)\right)
\]
is bijective. Let $K$ and $C$ denote the kernel and cokernel of $\Psi$, respectively. By Lemma~\ref{generator}, one has $K=0$. Since $V$ is projective, the abelian group $\Hom_{\mathcal{C}}(V,\,C)$ is the cokernel of $\Psi_*(V)$, hence vanishes. Applying Lemma~\ref{generator} again, we conclude that $C=0$. Therefore, $\Psi$ is an isomorphism.

We now obtain a chain of bijections
\begin{multline*}
 \varinjlim \Hom_\mathcal{C}(M \otimes N ,\, X_\alpha ) \simeq \varinjlim \Hom_\mathcal{C}(M  ,\, \Map(N,\, X_\alpha) ) \\
      \simeq\Hom_\mathcal{C}(M  ,\,  \varinjlim \Map(N,\, X_\alpha) ) \\
      \simeq\Hom_\mathcal{C}(M  ,\,  \Map(N,\, X) )
    \simeq\Hom_\mathcal{C}(M \otimes N ,\,  X )
\end{multline*}
showing that the induced homomorphism
\[
 \Psi_*(M \otimes N):    \varinjlim \Hom_\mathcal{C}(M \otimes N ,\, X_\alpha ) \to \Hom_\mathcal{C}(M \otimes N  ,\,X)
\]
is bijective. Hence $M \otimes N$ is finitely presented.

The same argument shows that if both $M$ and $N$ are finitely generated, then for any filtered system $\{X_\alpha\}$ whose transition morphisms are monomorphisms, the map $\Psi_*(M \otimes N)$ is bijective. Hence $M \otimes N$ is finitely generated. \qed
\end{proof}

\section{Nilpotent approximation of Smith ideals}
\label{sec:Smith}
\subsection{Monoidal structures of the arrow category of pointed symmetric monoidal model categories}
Let $[1]$ denote the category with two objects $0$ and $1$, and a single non-identity morphism $0 \to 1$. For any category $\mathcal{C}$, the {\it arrow category} $\mathrm{Ar}(\mathcal{C})$ of $\mathcal{C}$ is the category of covariant functors $[1] \to \mathcal{C}$. The category $\mathrm{Ar}(\mathcal{C})$ has evaluation functors $\mathrm{Ev}_{i}(f: X_0 \to X_1)=X_i \ (i=0,\,1)$, each admitting a left adjoint $L_i$ and a right adjoint $U_i$:
\[
\xymatrix@1{
\mathrm{Ev}_{i} : \mathrm{Ar}(\mathcal{C})  \ar[r]  & 
 \ar[l]<-1mm> \ar[l]<1mm>   \mathcal{C}:L_i ,\,U_i \ (i=0,\,1).
}
\]
If $\mathcal{C}$ is a pointed symmetric monoidal category, then $\mathrm{Ar}(\mathcal{C})$ carries two natural symmetric monoidal structures. Hovey~\cite[Theorem 1.2]{Smith-ideals} showed that these structures are induced by the monoidal structure on $\mathcal{C}$. For arrows $f:X_0 \to X_1$ and $g:Y_0 \to Y_1$, consider the commutative square
\[
 \xymatrix@1{ X_0 \otimes Y_0 \ar[r]^{f \otimes \mathrm{id}} 
\ar[d]_{\mathrm{id}\otimes g} & X_1 \otimes Y_0 
\ar[d]^{\mathrm{id}\otimes g}\\
X_0 \otimes Y_1 \ar[r]_{f \otimes \mathrm{id}}& X_1 \otimes Y_1.
}
\]
One of these structures is the {\it tensor product monoidal structure}, defined by
\[
f \otimes g = (f \otimes \mathrm{id}) \circ (
\mathrm{id} \otimes g)= ( \mathrm{id} \otimes g)\circ (f \otimes
\mathrm{id}) : X_0 \otimes Y_0 \to X_1 \otimes Y_1.
\]
The other is the {\it pushout product monoidal structure}, defined by the induced morphism
\[
 f \Box g: ( X_0 \otimes Y_1   ) \amalg_{ X_0 \otimes Y_0 } (   X_1 \otimes Y_0 ) \to  X_1 \otimes Y_1.
\]
Furthermore, the cokernel functor
\[
 \mathrm{cok}: \mathrm{Ar}(\mathcal{C}) \ni  (f :X \to Y) \mapsto (Y \to \mathrm{Coker}(f))  \in \mathrm{Ar}(\mathcal{C})
\]
is symmetric monoidal from the pushout product monoidal structure to the diagonal monoidal structure, and it admits a lax monoidal right adjoint
\[
 \mathrm{ker} : \mathrm{Ar}(\mathcal{C}) \ni  (f :X \to Y) \mapsto ( \mathrm{Ker}(f) \to X  )  \in \mathrm{Ar}(\mathcal{C})
\]
(see Hovey~\cite[Theorem 1.4]{Smith-ideals}). The right adjoint $\mathrm{ker}$ is called the {\it kernel} functor.

In the model-categorical setting, White--Yau~\cite{WY2024} develop the homotopy theory of commutative Smith ideals, namely commutative monoid objects in the pushout-product arrow category of a symmetric monoidal model category. Our notation $\CAlg(\mathrm{Ar}^\Box(\mathcal{M}))$ in Section~\ref{sec:model} follows this viewpoint.

\begin{definition}
A {\it Smith ideal} is a commutative monoid object in the arrow category $\mathrm{Ar}(\mathcal{C})$ equipped with the pushout product monoidal structure.
\end{definition}

\subsection{Nilpotent approximation of Smith ideals}
Let $\CAlg(\mathrm{Ar}^\Box(\mathcal{C}))$ denote the full subcategory of $\mathrm{Ar}(\mathcal{C})$ spanned by Smith ideals. We say that a Smith ideal $j:I \to A$ is nilpotent of degree $m\ge 1$ if the cokernel of the multiplication $\mu_m:I^{\otimes m} \to I$ is zero. This yields a filtration on $\CAlg(\mathrm{Ar}^\Box(\mathcal{C}))$: for each $n \ge 0$, let
$\CAlg^{n}(\mathrm{Ar}^\Box(\mathcal{C}))$ denote the full subcategory of $\CAlg(\mathrm{Ar}^\Box(\mathcal{C}))$ spanned by nilpotent Smith ideals of degree $n+1$. Then the inclusion functor
$\CAlg^{n}(\mathrm{Ar}^\Box(\mathcal{C}))
\to\CAlg(\mathrm{Ar}^\Box(\mathcal{C}))$ has a left adjoint given by
\[
\mathcal{P}^n:  \CAlg(\mathrm{Ar}^\Box(\mathcal{C})) \ni (j: I \to A) \mapsto (\tilde{j} : I/I^{n+1} \to A/I^{n+1}) \in \CAlg^{n}(\mathrm{Ar}^\Box(\mathcal{C})).
\]
Thus we obtain a tower of left adjoint functors
\[
  \cdots  \overset{\mathcal{P}^{n}}{\to}  \CAlg^{n}(\mathrm{Ar}^\Box(\mathcal{C})) \overset{\mathcal{P}^{n-1}}{\to} \CAlg^{n-1}(\mathrm{Ar}^\Box(\mathcal{C})) \overset{\mathcal{P}^{n-2}}{\to} \cdots \overset{\mathcal{P}^{0}}{\to} \CAlg^{0}(\mathrm{Ar}^\Box(\mathcal{C})),
\]
and the inverse-limit functor
\[
    \varprojlim{\mathcal{P}^n}: \CAlg(\mathrm{Ar}^\Box(\mathcal{C})) \to 2-\varprojlim\CAlg^{n}(\mathrm{Ar}^\Box(\mathcal{C})).
\]

We call the family of functors $\{\mathcal{P}^n\}$ the {\it nilpotent Taylor tower} for Smith ideals, the tower $\{\mathcal{P}^n(j) \}_{ n \ge 0}$ the {\it nilpotent Taylor tower} of $j$, and $\Lambda(j)= \varprojlim \mathcal{P}^n(j)$ the {\it nilpotent approximation} of $j$.

\begin{definition}
\label{analytic}
Let $\varphi: j \to j'$ be a morphism of Smith ideals. We say that $\varphi$ is an analytic equivalence if $\mathcal{P}^n(\varphi):\mathcal{P}^n(j) \to \mathcal{P}^n(j')$ is an isomorphism for each $n \ge 0$.
\end{definition}

\begin{proposition}
\label{analytic-limit-criterion}
Let $\varphi: j \to j'$ be a morphism of Smith ideals. Then $\varphi$ is an analytic equivalence if and only if the induced morphism
\[
\varprojlim \mathcal{P}^n(\varphi):\varprojlim \mathcal{P}^n(j) \to \varprojlim \mathcal{P}^n(j')
\]
is an isomorphism.
\end{proposition}
\begin{proof}
If $\varphi$ is an analytic equivalence, then the induced morphism on inverse limits is clearly an isomorphism. Conversely, assume that
\[
\varprojlim \mathcal{P}^n(\varphi):\varprojlim \mathcal{P}^n(j) \to \varprojlim \mathcal{P}^n(j')
\]
is an isomorphism. Then for each $m \ge 0$, the induced morphism
\[
\mathcal{P}^m(\varprojlim \mathcal{P}^n(\varphi)): \mathcal{P}^m(\varprojlim \mathcal{P}^n(j)) \to \mathcal{P}^m(\varprojlim \mathcal{P}^n(j'))
\]
is an isomorphism. Therefore the retraction $\mathcal{P}^m(\varphi):\mathcal{P}^m(j) \to \mathcal{P}^m(j')$ is also an isomorphism. Hence $\varphi$ is an analytic equivalence. \qed
\end{proof}

For any object $M$ of $\mathcal{C}$, write
$\mathcal{P}_{j}^n(M)=\mathcal{P}^n(j)\Box L_1(M): I/{I}^{n+1} \otimes M
\to A/I^{n+1} \otimes M$,
and call the inverse limit $\varprojlim \mathcal{P}_{j}^n(M)$ the {\it $j$-adic nilpotent approximation} of $M$.
\begin{definition}
\label{complete}
Let $j: I \to A$ be a Smith ideal. We say that $j$ is {\it complete} if $j \to \Lambda(j)$ is an analytic equivalence. An object $M$ is {\it $j$-adic complete} if $j \Box L_1(M) \to \varprojlim \mathcal{P}_j^n(M)$ is an analytic equivalence.
\end{definition}

\begin{definition}
\label{completion}
Let $\mathcal{C}$ be a pointed symmetric monoidal category. For any Smith ideal $j:I \to A$, let $\Lambda^\infty(j)$ denote the sequential colimit
\[
 \Lambda^\infty(j)=\varinjlim\left(\Lambda(j)\to\Lambda^2(j)\to\Lambda^3(j)\to\cdots\right).
\]
We call $\Lambda^\infty(j)$ the {\it completion} of $j$.
\end{definition}
Clearly, $\Lambda^\infty(j)$ is complete in the sense of Definition~\ref{complete}.

%
\begin{lemma}
\label{Split-End-c}
Let $\mathcal{C}$ be a locally presentable category, let $F: \mathcal{C} \to \mathcal{C}$ be an endofunctor, and let $\varphi: \mathrm{id}_\mathcal{C} \to F$ be a natural transformation.
Further, let $F^\infty$ denote the filtered colimit of the sequence of natural transformations
\[
  F \overset{F(\varphi)}{\to} F\circ F \overset{F^2(\varphi)}{\to}  F \circ F \circ F \overset{F^3(\varphi)}{\to} \cdots .
\]
Assume that the natural transformation $F(\varphi): F(-) \to F(F(-))$ is a split monomorphism. Then the canonical morphism $F(X) \to F^\infty(X)$ is an isomorphism for every finitely generated object $X$.
\end{lemma}
\begin{proof}
Let $\mathcal{C}[F^{-1}]$ be the localization defined as the $2$-colimit of the diagram
\[
   \mathcal{C} \overset{F}{\to} \mathcal{C}\overset{F}{\to} \mathcal{C}\overset{F}{\to} \mathcal{C} \overset{F}{\to} \cdots\,.
\]
Then $F^\infty$ is left adjoint to the inclusion functor $\mathcal{C}[F^{-1}] \to \mathcal{C}$, and hence preserves all small colimits. If $X$ is finitely generated, then $F^\infty(X)$ is also finitely generated. Indeed, for any filtered system $\{M_\alpha\}$ whose transition morphisms are monomorphisms, the map
\[
\varinjlim \Hom(F^{\infty}(X),\,M_\alpha ) \to \Hom_\mathcal{C}(F^\infty(X),\, \varinjlim M_\alpha )
\]
is a retract of the composite isomorphism
\begin{multline*}
\varinjlim \Hom_\mathcal{C}(X,\,F^\infty(M_\alpha) ) \simeq \varinjlim \Hom_{\mathcal{C}[F^{-1}]}(F^{\infty}(X),\, F^\infty( M_\alpha) ) \\
\to \Hom_{\mathcal{C}[F^{-1}]}(F^\infty(X),\, \varinjlim F^\infty(M_\alpha) ) \simeq \Hom_\mathcal{C}(X,\, \varinjlim F^\infty(M_\alpha) ) .
\end{multline*}
Since $F^\infty(X)$ is a filtered colimit of split monomorphisms $F^n(X)
\to F^{n+1}(X) \ ( n \ge 1)$, Proposition~\ref{f-g} implies that the identity on
$F^\infty(X)$ factors through $F^{n}(X)$ for some $n$. Therefore the monomorphism $F^n(\varphi)(X): F^n(X) \to F^{n+1}(X)$ is an isomorphism, and hence so is its retract $F(\varphi)(X): F(X) \to F(F(X))$. It follows that the canonical map $F(X) \to F^\infty(X)$ is an isomorphism. \qed
\end{proof}

\begin{proposition}
Let $j:I \to A$ be a finitely generated Smith ideal. Then $\Lambda(j)$ is already complete.
\end{proposition}
\begin{proof}
Let $\varphi: j \to \varprojlim \mathcal{P}^n (j) = \Lambda(j)$ denote the canonical morphism induced by the localization maps $\varphi^n(j):j \to \mathcal{P}^n(j)$ for each $n \ge 1$. Since $\mathcal{P}^n(j) \to \mathcal{P}^n(\mathcal{P}^n(j))$ is an isomorphism for every $n \ge 1$, the identity of $\mathcal{P}^n(j)$ factors through $\mathcal{P}^n (\varprojlim \mathcal{P}^n (j))$. Passing to inverse limits, we obtain a split monomorphism $\Lambda(\varphi): \Lambda(j) \to \Lambda(\Lambda(j))$. Applying Lemma~\ref{Split-End-c} to the finitely generated Smith ideal $j$, we conclude that $\Lambda(\varphi)$ is an isomorphism. Hence $\Lambda(j)$ is complete. \qed
\end{proof}

Finitely generated (resp. finitely presented) objects of the arrow category $\mathrm{Ar}(\mathcal{C})$ of a locally presentable category $\mathcal{C}$ are defined termwise; see Ad{\'a}mek and Rosick{\'y}~\cite[p.44, Example (1)]{LocPresentable}. Concretely, let $f: C \to D$ be an arrow between finitely generated objects, and let $\{ g_\alpha: X_\alpha \to Y_\alpha \}$ be a filtered inductive system of arrows whose transition morphisms $X_\alpha \to X_\beta$ and $Y_\alpha \to Y_\beta \ (\beta \ge \alpha)$ are monomorphisms. Write $X= \varinjlim X_\alpha$ and $Y= \varinjlim Y_\alpha$. Since finite limits commute with filtered colimits in $\mathcal{C}$, the canonical map
\[
 \varinjlim \left(\Hom_\mathcal{C}(C,\,X_\alpha) \times_{\Hom_{\mathcal{C}}(C,\,Y_\alpha)} \Hom_\mathcal{C} (D,\,Y_\alpha)  \right) \to  \Hom_\mathcal{C}(C,\,X) \times_{\Hom_{\mathcal{C}}(C,\,Y)} \Hom_\mathcal{C} (D,\,Y)
\]
is bijective. Therefore every morphism $s: f \to g$ factors through some $g_\alpha: X_\alpha \to Y_\alpha$. The same argument, without assuming that the transition morphisms are monomorphisms, shows that if $C$ and $D$ are finitely presented, then $f: C \to D$ is finitely presented.

\begin{proposition}
\label{fg-push} Let $\mathcal{C}$ be a locally presentable symmetric
monoidal abelian category generated by its monoidal unit $V$. Assume that
$V$ is finitely presented and projective. Then the class of finitely
generated (resp. finitely presented) objects of the arrow category
$\mathrm{Ar}(\mathcal{C})$ is closed under pushout products.
\end{proposition}
\begin{proof}
Let $f: X_0 \to X_1$ and $g: Y_0 \to Y_1$ be finitely generated (resp. finitely presented) arrows. By Proposition~\ref{closed-compact}, each tensor product $X_i \otimes Y_j$ for $0 \le i,\,j \le 1$ is a finitely generated (resp. finitely presented) object of $\mathcal{C}$. Since filtered colimits commute with finite limits, finitely generated (resp. finitely presented) objects of $\mathcal{C}$ are closed under finite colimits. Therefore the pushout $X_0 \otimes Y_1 \amalg_{X_0 \otimes Y_0} X_1 \otimes Y_0$ is finitely generated (resp. finitely presented). Hence the pushout product
\[
 f \Box g: X_0 \otimes Y_1 \amalg_{X_0 \otimes Y_0} X_1 \otimes Y_0 \to X_1 \otimes Y_1
\]
is finitely generated (resp. finitely presented). \qed
\end{proof}

\begin{lemma}
\label{compact-comp-c}
Let $j: I \to A$ be a finitely generated Smith ideal of $A$, and let $M$ be a finitely generated object of the symmetric monoidal category $\mathcal{C}$. Assume that the monoidal unit $V$ generates $\mathcal{C}$. Then $\Lambda_j(M)$ is $j$-adic complete.
\end{lemma}
\begin{proof}
Since finitely generated objects of any compactly generated closed symmetric monoidal category generated by the monoidal unit are closed under tensor products, Proposition~\ref{fg-push} implies that the pushout product $j\Box L_0(M)$ is finitely generated. By Lemma~\ref{Split-End-c}, therefore, the canonical map $\Lambda_j(M) \to \varinjlim_n \Lambda^n_j(M)$ is an isomorphism. Hence the $j$-adic nilpotent approximation $\Lambda_j(M)$ is complete. \qed
\end{proof}

\begin{theorem}
\label{TheoremA1}
Let $\mathcal{C}$ be a locally presentable closed symmetric monoidal abelian category such that the monoidal unit $V$ is a finitely generated projective generator. Let $j: I \to A$ be a finitely generated Smith ideal, and let $M$ be a $\kappa$-filtered colimit of finitely generated objects of $\mathcal{C}$. Then the $j$-adic nilpotent approximation $\Lambda_j(M)$ is $j$-adic complete.
\end{theorem}
\begin{proof}
Since $\kappa$-filtered colimits commute with $\kappa$-small limits in $\mathcal{C}$ (see Ad{\'a}mek and Rosick{\'y}~\cite[p.44, Proposition 1.59]{LocPresentable}), the theorem reduces to the case where $M$ is finitely generated. The claim in that case follows from Lemma~\ref{compact-comp-c}. \qed
\end{proof}

Let $\mathcal{P}_j^n(\mathcal{C})$ denote the full replete subcategory of $\mathrm{Ar}(\mathcal{C})$ given by the essential image of the functor $\mathcal{P}_j^n:  \mathcal{C} \to  \mathrm{Ar}(\mathcal{C})$ for each $n \ge 0$, and let $\Lambda_j(\mathcal{C})$ be the full replete subcategory given by the essential image of the $j$-adic nilpotent approximation functor $\Lambda_j :\mathcal{C} \to \mathrm{Ar}(\mathcal{C})$.
For each $n \ge 1$, let
\[
\tau_n:\mathcal{P}_j^n(\mathcal{C})\to \mathcal{P}_j^{n-1}(\mathcal{C})
\]
be the restriction of the canonical truncation morphism $\mathcal{P}_j^n(M)\to \mathcal{P}_j^{n-1}(M)$. We write $2-\varprojlim\mathcal{P}_j^n(\mathcal{C})$ for the $2$-limit of the tower
\[
\cdots\xrightarrow{\tau_{n+1}}\mathcal{P}_j^n(\mathcal{C})\xrightarrow{\tau_n}\mathcal{P}_j^{n-1}(\mathcal{C})\xrightarrow{}\cdots\xrightarrow{}\mathcal{P}_j^0(\mathcal{C}).
\]
Then we have the following:
\begin{theorem}
\label{TheoremA2}
Let $\mathcal{C}$ be a locally presentable closed symmetric monoidal abelian category.  Assume that the monoidal unit is a finitely generated projective generator.  Let $j: I \to A$ be a finitely generated Smith ideal. 
Then the functor  
\[
  \varprojlim : 2-\varprojlim  \mathcal{P}_j^n(\mathcal{C})     \to \Lambda_j(\mathcal{C} ) 
\]  
is a categorical equivalence.
\end{theorem}

\begin{proof}
Define
\[
\Phi: \mathcal{C}\to 2-\varprojlim\mathcal{P}_j^n(\mathcal{C}),\qquad M\mapsto \{\mathcal{P}_j^n(M)\}_{n\ge 0},
\]
where the structure maps are given by the canonical truncations. By construction,
\[
(\varprojlim\circ \Phi)(M)=\varprojlim_n\mathcal{P}_j^n(M)=\Lambda_j(M).
\]
Hence $\varprojlim$ is essentially surjective onto $\Lambda_j(\mathcal{C})$ by definition of essential image.

To prove full faithfulness, let $\mathbf{X}=\{X_n\}$ and $\mathbf{Y}=\{Y_n\}$ be objects of $2-\varprojlim\mathcal{P}_j^n(\mathcal{C})$. Since each $\mathcal{P}_j^n(\mathcal{C})$ is an essential image, we can choose $M,N\in\mathcal{C}$ and compatible isomorphisms
\[
X_n\simeq \mathcal{P}_j^n(M),\qquad Y_n\simeq \mathcal{P}_j^n(N) \quad (n\ge 0).
\]
By Theorem~\ref{TheoremA1}, the canonical morphisms
\[
j\Box L_1(M)\to \Lambda_j(M),\qquad j\Box L_1(N)\to \Lambda_j(N)
\]
are analytic equivalences. Applying Proposition~\ref{analytic-limit-criterion} to these maps, for each $n\ge 0$ we obtain canonical identifications
\[
\mathcal{P}_j^n(M)\simeq \mathcal{P}_j^n(\Lambda_j(M)),\qquad
\mathcal{P}_j^n(N)\simeq \mathcal{P}_j^n(\Lambda_j(N)).
\]
Therefore a morphism $\Lambda_j(M)\to\Lambda_j(N)$ determines and is determined by a compatible family
\[
\{\mathcal{P}_j^n(M)\to \mathcal{P}_j^n(N)\}_{n\ge 0},
\]
that is, by a morphism $\mathbf{X}\to\mathbf{Y}$ in $2-\varprojlim\mathcal{P}_j^n(\mathcal{C})$. Hence $\varprojlim$ is fully faithful.

Hence, $\varprojlim$ is a categorical equivalence. \qed
\end{proof}

\subsection{Almost nilpotent approximation of almost modules}
\label{sec:almost}
We briefly recall the basic notions of almost mathematics. Throughout this subsection, we fix a unital commutative ring $V$ and an idempotent ideal $\mathfrak{m}$ of $V$. We always assume that $\tm:= \mathfrak{m} \otimes_V
\mathfrak{m}$ is a flat $V$-module.

A $V$-module $M$ is said to be {\it almost zero} if $\mathfrak{m}M=0$. The full subcategory of almost zero $V$-modules is a Serre subcategory of the category of $V$-modules, and we let $\AMod_V$ denote the corresponding localization.

A $V$-module $M$ is said to be {\it firm} if the product $\mu_M: \mathfrak{m} \otimes_V M \to M$ is an isomorphism. It is well known that, for every $V$-module $M$, the module $\tm \otimes_V M$ is firm. Let $\FMod_V$ denote the full subcategory of $\Mod_V$ spanned by firm $V$-modules. Quillen~\cite{Q} proved that the functor
\[
 \tm \otimes_V (-) : \Mod_V\to \FMod_V 
\]   
induces a categorical equivalence $\AMod_V \simeq \FMod_V$. Therefore a homomorphism $f: M \to N$ is an isomorphism in the category $\AMod_V$ if and only if the induced map $\tm \otimes_V f: \tm \otimes_V M \to \tm \otimes_V N$ is an isomorphism of $V$-modules. We say that a homomorphism $f: M \to N$ of $V$-modules is {\it almost injective} (resp. {\it almost surjective}, resp. {\it an almost isomorphism}) if the induced morphism $\tm \otimes_V f: \tm \otimes_V M \to \tm \otimes_V N$ is injective (resp. surjective, resp. an isomorphism).

Let $j: I \to A$ be an almost finitely generated Smith ideal, and let $M$ be an almost $A$-module. If $j \otimes_A M: I \otimes_A M \to M$ induces an almost isomorphism $IM \to M$, then $\mathcal{P}^n(j)\otimes_A M$ is almost zero for every $n \ge 0$.

\begin{definition}
Let $j: I \to A$ be an almost Smith ideal and $M$ an almost $A$-module. The {\it almost nilpotent approximation} of $M$ is defined to be the inverse limit of the morphisms
\[
  \mathcal{P}^n(j \otimes_A M): IM/ I^{n+1}M  \to M / I^{n+1} M
\]
as $n$ ranges over all nonnegative integers.
\end{definition}

The following proposition is a special case of Theorem~\ref{TheoremA1}:
\begin{proposition}
\label{almost-f-g}
Let $A$ be an almost $V$-algebra and $j: I \to A$ a Smith ideal such that $I$ is an almost finitely generated $A$-module. Then, for any almost $A$-module $M$, the almost nilpotent approximation $\Lambda_j(M)$ is almost complete.
\end{proposition}

\begin{proof}
This is a special case of Theorem~\ref{TheoremA1}. \qed
\end{proof}

\begin{corollary}[{\rm c.f.} \cite{compFG} Corollary 3.6]
Let $A$ be a commutative ring and $j: I \to A$ a Smith ideal such that $I$ is an almost finitely generated $A$-module. Then, for any almost $A$-module $M$, the almost nilpotent approximation $\Lambda_j(M)$ is complete.
\end{corollary}
\begin{proof}
This is the special case of Proposition~\ref{almost-f-g} with $V=\mathfrak{m}=A$. \qed
\end{proof}

\begin{definition}
Let $A$ be an almost $V$-algebra and $M$ an almost $A$-module. We say that $M$ is {\it almost compactly projective} if $\tm \otimes_{V} M$ is a direct factor of $\tm \otimes_V A^n$ for some $n \ge 0$. Let $\APMod_A^{\rm ac}$ denote the full subcategory of $\AMod_A$ spanned by almost compactly projective $A$-modules.
\end{definition}

\begin{theorem}[{\rm c.f.} \cite{GR} p.145, Lemma 5.3.16 and p.147, Theorem 5.3.24]
\label{MainThm}
Let $A$ be an almost $V$-algebra and $j: I \to A$ an almost finitely generated Smith ideal of $A$. Then the functor
\[
  \varprojlim : 2-\varprojlim \mathcal{P}^n_j(\AMod_A) \to \Lambda_{j}(\AMod_A)  
\]
is a categorical equivalence. Furthermore, restricting to almost compactly projective modules yields a categorical equivalence
\[
  \varprojlim : 2-\varprojlim \mathcal{P}^n_j(\APMod^{\rm ac}_A)   \to  \APMod^{\rm ac}_{\hat{j}},
\]
where $\APMod^{\rm ac}_{\hat{j}}$ denotes the essential image of $j \Box L_1(-): \APMod^{\rm ac}_{\hat{A}} \to \mathrm{Ar}(\APMod^{\rm ac}_{\hat{A}})$.
\end{theorem}
\begin{proof}
The first statement is a special case of Theorem~\ref{TheoremA2}. For the second statement, it is enough to show that every almost compactly projective firm $\hat{A}$-module is already $\hat{I}$-adic almost complete. By definition, such a module is a direct factor of $\tm \otimes_V \hat{A}^n$ for some $n \ge 0$. Therefore it is already almost complete. \qed
\end{proof}

\section{Nilpotent approximation of Smith ideals in symmetric monoidal model categories}
\label{sec:model}

\subsection{Homotopically nilpotent Smith ideals}

Let $\mathcal{M}$ be a stable symmetric monoidal model category. Throughout this section, every symmetric monoidal model category is assumed to be closed, proper, and combinatorial. By Hovey~\cite[Theorem 2.1 (4), Theorem 3.1 (5), and {\rm c.f.} Theorem 1.4]{Smith-ideals}, the arrow category $\mathrm{Ar}^\Box(\mathcal{M})$ with the injective model structure is symmetric monoidal with respect to the pushout product.
%

Smith ideals were defined in Section~\ref{sec:Smith}. A Smith ideal $j: I \to A$ is said to be homotopically nilpotent of degree $m \ge 1$ if the multiplication $\mu_m: I^{\otimes m} \to I$ is null-homotopic. Equivalently, the homotopy cokernel morphism $I \to \mathrm{Coker}(\mu_m)$ is a weak equivalence.

We now introduce a Taylor tower for Smith ideals. When the transferred model structure on commutative monoids in $\mathrm{Ar}^\Box(\mathcal{M})$ exists, we denote it by $\CAlg(\mathrm{Ar}^\Box(\mathcal{M}))$ and refer to it as the model category of Smith ideals; for the motivic symmetric spectra used later, this follows from White--Yau~\cite[Theorem~A and Example~4.4.2]{WY2024}. Since the Smith-ideal model category is proper and combinatorial in the cases considered here, the left Bousfield localization imposing the null-homotopy of the $(n+1)$-fold multiplication exists. Let $\CAlg^{n}(\mathrm{Ar}^\Box(\mathcal{M}))$ denote its full subcategory spanned by Smith ideals that are nilpotent of degree $n+1$. Then the inclusion functor $\CAlg^{n}(\mathrm{Ar}^\Box(\mathcal{M})) \to\CAlg(\mathrm{Ar}^\Box(\mathcal{M}))$ has a left adjoint
\[
\mathcal{P}^n:  \CAlg(\mathrm{Ar}^\Box(\mathcal{M})) \to  \CAlg^{n}(\mathrm{Ar}^\Box(\mathcal{M})), 
\]
which is the left Bousfield localization with respect to the class of morphisms $I \to
\mathrm{Coker}(\mu_{n+1}: I^{\otimes n+1} \to I )$. Let $I^{m}$ denote the homotopy image of $\mu_{m}: I^{\otimes m} \to I$, and let $I/I^{m}$ denote the homotopy cokernel of $I^{m} \to I$ for $m \ge 0$. By definition of the model category $\CAlg^{n}(\mathrm{Ar}^\Box(\mathcal{M}))$, the canonical morphism $I \to I/I^{n+1}$ induces a weak equivalence $\mathcal{P}^n(I) \to \mathcal{P}^n(I/I^{n+1})$. Furthermore, since $\mathcal{M}$ is stable, the homotopy cocartesian square
\[
 \xymatrix@1{
  I^m \ar[r] \ar[d] & I^m/I^{n+1} \ar[d] \\
  I \ar[r]  & I /I^{n+1}
}
\]
is also homotopy cartesian. Therefore the induced morphism $I^m \to I^m/I^{n+1}$ is a weak equivalence if and only if $I \to I/I^{n+1}$ is a weak equivalence. In particular, $\mathcal{P}^n(I^n) \to \mathcal{P}^n(I^n/I^{n+1})$ is a weak equivalence for each $n \ge 1$.
\begin{lemma}
\label{graded}
Let $\mathcal{M}$ be a stable symmetric monoidal model category and $j: I \to A$ a Smith ideal. For any $n \ge 1$, the left localization functor
\[
 \mathcal{P}^n: \CAlg(\mathrm{Ar}^\Box(\mathcal{M})) \to \CAlg^{n}(\mathrm{Ar}^\Box(\mathcal{M}))
\]
induces the following homotopy fiber sequence in the stable model category $\mathrm{Ar}^\Box(\mathcal{M})$:
\[
      L_{1} (I^n/I^{n+1}) \to \mathcal{P}^n(j) \to  \mathcal{P}^{n-1}(j).
\]
\end{lemma}
\begin{proof}
Indeed, one has the commutative diagram
\[
 \xymatrix@1{
  I^n/I^{n+1} \ar[r] \ar[d] & I/I^{n+1} \ar[d]  \ar[r] & I/I^n  \ar[d]\\
  I^n/I^{n+1} \ar[r]  & A /I^{n+1}  \ar[r] & A/I^n,
 }
\]
where both horizontal rows are homotopy fiber sequences. The vertical maps represent $\mathcal{P}^n(j)$ and $\mathcal{P}^{n-1}(j)$ in the localization $\CAlg^{n}(\mathrm{Ar}^\Box(\mathcal{M}))$. Hence the homotopy fiber of $\mathcal{P}^n(j) \to \mathcal{P}^{n-1}(j)$ is equivalent to $L_1(I^n/I^{n+1})$. \qed
\end{proof}
\begin{lemma}
\label{n/n+1}
Let $\mathcal{M}$ be a stable symmetric monoidal model category and $j: I \to A$ a Smith ideal in $\mathcal{M}$. Then $\mathrm{Coker}(j) \otimes I^n \to \mathrm{Coker}(I^{n+1} \to I^n)$ is a weak equivalence.
\end{lemma}
\begin{proof}
Fix an integer $n \ge 0$. There is a homotopy cokernel sequence $I \otimes I^n \to A \otimes I^n \to \mathrm{Coker}(j) \otimes I^n$. Since $\mathcal{M}$ is stable, the unit $u: \mathrm{Id}_\mathcal{M} \to \mathrm{ker} \circ \mathrm{cok}$ of the Quillen adjunction $\mathrm{cok} :\mathrm{Ar}(\mathcal{M}) \rightleftarrows \mathrm{Ar}(\mathcal{M}) : \mathrm{ker}$ is a weak equivalence by Hovey's theorem~\cite[Theorem 4.3]{Smith-ideals}. Therefore the induced morphisms $I^n \otimes I \to I^{n+1}$ and $I^n \otimes A \to I^n$ are weak equivalences, and hence so is the induced morphism $\mathrm{Coker}(j) \otimes I^n \to \mathrm{Coker}(I^{n+1} \to I^n)$. \qed
\end{proof}

We call the homotopy limit of functors
\[
\varprojlim \mathcal{P}^n:
\CAlg(\mathrm{Ar}^\Box(\mathcal{M})) \to 2-\varprojlim
\CAlg^{n}(\mathrm{Ar}^\Box(\mathcal{M}))
\]
the {\it homotopically nilpotent approximation} of Smith ideals.

\begin{definition}
\label{m-analytic}
Let $\varphi: j \to j'$ be a morphism of Smith ideals. We say that $\varphi$ is a {\it homotopically analytic equivalence} if $\mathcal{P}^n(\varphi):\mathcal{P}^n(j) \to \mathcal{P}^n(j')$ is a weak equivalence for each $n \ge 0$.
\end{definition}

\begin{proposition}
Let $\varphi: j \to j'$ be a morphism of Smith ideals. Then $\varphi$ is a homotopically analytic equivalence if and only if the induced morphism
\[
\varprojlim\mathcal{P}^n(\varphi):\varprojlim \mathcal{P}^n(j) \to \varprojlim \mathcal{P}^n(j')
\]
is a weak equivalence.
\end{proposition}  
\begin{proof}
If $\varphi$ is a homotopically analytic equivalence, then the induced morphism on homotopy limits is a weak equivalence because homotopy limits preserve levelwise weak equivalences. Conversely, assume that
\[
\varprojlim\mathcal{P}^n(\varphi):\varprojlim \mathcal{P}^n(j) \to \varprojlim \mathcal{P}^n(j')
\]
is a weak equivalence. Then for each $m \ge 0$, the induced morphism
\[
\mathcal{P}^m(\varprojlim(\mathcal{P}^n(\varphi))): \mathcal{P}^m(\varprojlim \mathcal{P}^n(j)) \to \mathcal{P}^m(\varprojlim \mathcal{P}^n(j'))
\]
is again a weak equivalence. Hence the homotopy retraction $\mathcal{P}^m(\varphi):\mathcal{P}^m(j) \to \mathcal{P}^m(j')$ is a weak equivalence. Therefore $\varphi$ is a homotopically analytic equivalence. \qed
\end{proof}

\begin{definition}
\label{m-w-compact}
An object $X$ of $\mathcal{M}$ is said to be {\it compact} if, for any filtered inductive system $\{ M_\alpha \}$, the induced map
\[
  \varinjlim \Hom_\mathcal{M}(X,\,  M_{\alpha} ) \to   \Hom_\mathcal{M}(X,\, \varinjlim \ M_{\alpha} )
\]
is an isomorphism in the homotopy category $\mathrm{Ho}(\mathcal{M})$. We say that $X$ is {\it weakly compact} if, for any filtered inductive system $\{ M_\alpha \}$ whose transition morphisms are homotopically splitting monomorphisms, the same induced map is an isomorphism in $\mathrm{Ho}(\mathcal{M})$.
\end{definition}

Let $\mathcal{M}$ be a closed left proper combinatorial model category. A morphism $f:X \to Y$ is said to be a {\it homotopical epimorphism} if, for every object $A$, the induced map
\[
  f^*:  \Hom_{\mathrm{Ho}(\mathcal{M})}(Y,\,A ) \to \Hom_{\mathrm{Ho}(\mathcal{M})}(X,\,A )
\]
is injective, where $\mathrm{Ho}(\mathcal{M})$ denotes the homotopy category of $\mathcal{M}$.
\begin{lemma}
\label{w-compact}
Let $\mathcal{M}$ be a left proper combinatorial model category. Let $X$ be an object. If $X$ admits a homotopical epimorphism $f: C \to X$ from a compact object $C$, then $X$ is weakly compact.
\end{lemma}
\begin{proof}
Let $\{M_\alpha,\,i_{\beta \alpha}\}$ be a filtered inductive system whose transition morphisms are homotopically splitting monomorphisms, and let $M$ denote the homotopy colimit. Clearly, the induced map
\[
  \varinjlim \Hom_{\mathrm{Ho}(\mathcal{M})}(X,\,  M_{\alpha} ) \to  \Hom_{\mathrm{Ho}(\mathcal{M})}  (X,\, M )
\]
is injective. Given a morphism $\psi: X \to M$, there exists $\alpha$ such that $\psi \circ f:C \to M$ factors through some $M_\alpha$. Let $i_\alpha: M_\alpha \to M$ be the canonical morphism. Then $i_\alpha$ has a homotopy retraction $p_\alpha: M \to M_\alpha$. Since $f$ is a homotopical epimorphism, the composition $i_\alpha \circ (p_\alpha \circ \psi): X \to M_\alpha \to M$ is homotopic to $\psi$. Therefore $\psi$ is represented by a morphism factoring through $M_\alpha$, and the above induced map is surjective. \qed
\end{proof}

\begin{lemma}
\label{Split-End}
Let $\mathcal{M}$ be a combinatorial model category, and let $F: \mathcal{M} \to \mathcal{M}$ be an endofunctor preserving all small colimits together with a natural transformation $\varphi: \mathrm{id}_\mathcal{M} \to F$.
Further, let $F^\infty$ denote the filtered colimit of the natural transformations
\[
  F \overset{F(\varphi)}{\to} F\circ F \overset{F^2(\varphi)}{\to}  F \circ F \circ F \overset{F^3(\varphi)}{\to} \cdots.
\]
Assume that for every compact object $X$, the morphism $F(\varphi): F(X) \to F(F(X))$ is a homotopically splitting monomorphism. Then the canonical natural transformation $F(X) \to F^\infty(X)$ is a weak equivalence for every weakly compact object $X$.
\end{lemma}
\begin{proof}
Since $F^\infty$ is left adjoint to the inclusion functor
$\mathcal{M}[F^{-1}] \to \mathcal{M}$, if $X$ is weakly compact,
then $F^\infty(X)$ is also weakly compact. Moreover,
$F^\infty(X)$ is a filtered colimit of splitting monomorphisms $F^n(X)
\to F^{n+1}(X) \ ( n \ge 1)$. Hence the identity on
$F^\infty(X)$ factors through $F^{n}(X)$ for some $n$. Therefore the homotopically splitting monomorphism $F^n(\varphi)(X): F^n(X) \to F^{n+1}(X)$ is a weak equivalence. It follows that $F(\varphi)(X): F(X) \to F(F(X))$, being a retract of $F^n(\varphi)(X)$, is also a weak equivalence. \qed
\end{proof}

\begin{proposition}
\label{compact-comp}
Let $j:I \to A$ be a Smith ideal, and write $\Lambda(j)= \varprojlim \mathcal{P}^n (j)$. If $j$ is weakly compact, then $\Lambda(j)$ is homotopically complete.
\end{proposition}
\begin{proof}
Let $\varphi: j \to \varprojlim \mathcal{P}^n (j) = \Lambda(j)$ denote the
canonical morphism induced by the localization $\varphi^n(j):j \to
\mathcal{P}^n(j)$ for each $n \ge 1$. Since $\mathcal{P}^n(j) \to \mathcal{P}^n(\mathcal{P}^n(j))$ is an
isomorphism for every $n \ge 1$, the identity of $\mathcal{P}^n(j)$ factors
through $\mathcal{P}^n (\varprojlim \mathcal{P}^n (j))$. Passing to filtered limits, we
obtain a homotopically splitting monomorphism $\Lambda(\varphi): \Lambda(j) \to
\Lambda(\Lambda(j))$ for any Smith ideal $j$. Applying
Lemma~\ref{Split-End} to the weakly compact Smith ideal $j$, we conclude
that $\Lambda(\varphi): \Lambda(j) \to \Lambda(\Lambda(j))$ is a weak equivalence. Hence $\Lambda(j)$ is homotopically complete. \qed
\end{proof}

\begin{remark}
The logical flow of Section~\ref{sec:model} is as follows: Lemma~\ref{graded} identifies the graded homotopy fiber sequence for the nilpotent tower, Lemma~\ref{n/n+1} compares graded pieces with $\mathrm{Coker}(j)$-modules, Definition~\ref{m-w-compact} and Lemma~\ref{w-compact} provide the weak-compactness criterion, and Lemma~\ref{Split-End} supplies the formal convergence step. Combining these ingredients yields Proposition~\ref{compact-comp}.
\end{remark}

\section{Application to algebraic cobordism}
\label{sec:MSP}
In this section, we first recall motivic spaces, motivic spectra, and motivic $\E_\infty$-rings. We fix a base scheme $S$, and let $\Sm_S$ denote the small category of smooth $S$-schemes.

\subsection{Motivic spaces and motivic spectra}

Let $\mathcal{S}_S:=\SShv^{\A^1}(\Sm_S)_*^{\rm Nis}$ denote the simplicial model category of pointed simplicial sheaves on the Nisnevich site $(\Sm_S)_{\rm Nis}$, equipped with the $\A^1$-model structure introduced by Morel--Voevodsky~\cite[p.106, Section 3.2, Definition 2.1]{MV}. An object of $\mathcal{S}_S$ is called a {\it pointed motivic space}. The $\A^1$-model structure on $\mathcal{S}_S$ is proper, closed simplicial, and combinatorial. The smash product $-\wedge -$ on $\mathcal{S}_S$ naturally induces a symmetric monoidal model structure.

For any pointed motivic space $T$, let $\mathbf{Spt}^\Sigma_{T}(S)$ denote the stable model category of symmetric $T$-spectra of pointed motivic spaces in the sense of Hovey~\cite{HoveySpectra} and Jardine~\cite[Section 4.1]{Jardine}. This model structure is the stable model structure of Jardine~\cite[p.516, Theorem 4.15]{Jardine}; it is proper, combinatorial, and closed simplicial. A {\it motivic spectrum} is an object of $\mathbf{Spt}^\Sigma_{\mathbb{P}^1_+}(S)$, where $\mathbb{P}^1_+$ denotes the pointed projective line. Moreover, by Jardine~\cite[pp.520--523, Propositions 4.19 and 4.20]{Jardine} (see also \cite[Theorem A.38]{MR2597741}), this stable model category is also a proper combinatorial closed symmetric monoidal model category with respect to the smash product.

In this section, for any motivic spectrum $M$, we write
\[
 M[n]= \begin{cases} 
(\mathbb{P}^1_+)^{\wedge n}\wedge M \quad (n \ge 0),\\
\Omega_{\mathbb{P}^1}^{-n}M \quad (n < 0),
       \end{cases}
\]
where $\Omega_{\mathbb{P}^1}^{m}M$ denotes the $m$-fold $\mathbb{P}^1$-loop object, namely the derived internal Hom object written as $\Map((\mathbb{P}^1_+)^{\wedge m},\,M)$, for $m>0$.
Then we have a Quillen adjunction
\[
  \Sigma_+^\infty: \SShv^{\A^1}(\Sm_S)^{\rm Nis}_* \rightleftarrows \mathbf{Spt}^\Sigma_{\mathbb{P}^1_+}(S) : \Omega_+^\infty.
\]
Here, the right Quillen functor $\Omega_+^\infty:\mathbf{Spt}^\Sigma_{\mathbb{P}^1_+}(S) \to \SShv^{\A^1}(\Sm_S)^{\rm Nis}_*$ is the evaluation-at-level-$0$ functor. That is, for any motivic spectrum $M=\left(M_0,\,M_1,\,\ldots\right)$, the functor $\Omega_+^\infty$ sends $M$ to $M_0$. The left Quillen functor $\Sigma_+^\infty : \SShv^{\A^1}(\Sm_S)^{\rm Nis}_* \to \mathbf{Spt}^\Sigma_{\mathbb{P}^1_+}(S)$ is the usual $\mathbb{P}^1_+$-suspension spectrum functor.

\begin{definition}
A {\it motivic $\E_\infty$-ring} is a homotopically commutative monoid object in the stable symmetric monoidal model category of motivic spectra.
\end{definition}

We also recall the notion of periodicity for motivic spectra.
\begin{definition}
We say that a motivic $\E_\infty$-ring $R$ is {\it periodic} if there exists an element $\beta \in \pi_0(\Map(\mathbb{P}^1_+,\,R))$ such that $\beta$ is a unit in the graded ring $\bigoplus_{n \in \Z} \pi_0(R[n])$.
\end{definition}

\begin{proposition}
\label{periodic-app}
Let $j: I \to A$ be the homotopy kernel of a morphism $f: A \to R$ of motivic $\E_\infty$-rings. If $R$ is periodic, then for each $n \ge 0$, the nilpotent approximation $\mathcal{P}^{n}(j)$ of degree $n+1$ and the $j$-adic nilpotent approximation $\Lambda_{j}(A)$ are also periodic.
\end{proposition}
\begin{proof}
We prove the proposition by induction on $n \ge 0$. The case $n=0$ is clear. Assume that $\mathcal{P}^{n-1}(j)$ is periodic.
By Lemma~\ref{graded}, the homotopy fiber of the canonical morphism
\[
\mathcal{P}^{n}(j)  \to \mathcal{P}^{n-1}(j)
\]
is weakly equivalent to $L_1(I^n/I^{n+1})$. By Lemma~\ref{n/n+1}, the morphism $R \otimes I^n \to I^n/I^{n+1}$ is a motivic equivalence, so the term $L_1(I^n/I^{n+1})$ is periodic for every $n \ge 0$. Hence $\mathcal{P}^{n}(j)$ is periodic. \qed
\end{proof}

For any motivic spectrum $M$, the homotopy coproduct
\[
  PM = \bigvee_{n \in \Z} M[n]
\] 
is periodic by construction and is called the {\it periodization} of $M$.
\begin{proposition}
\label{split-P}
Let $j: I \to A$ be a Smith ideal of a motivic $\E_\infty$-ring whose cokernel $R$ is a periodic motivic $\E_\infty$-algebra. Then, for any $n \ge 0$, the nilpotent approximation $\mathcal{P}^n(j)$ of degree $n+1$ is a retract of the periodization $Pj: PI \to PA$.
\end{proposition}
\begin{proof}
By the universal property of periodization, the canonical morphism
$M\to PM$ is initial among morphisms from $M$ to periodic motivic spectra;
equivalently, for every periodic motivic spectrum $N$, restriction along
$M\to PM$ induces a weak equivalence
$\Map(PM,\,N)\simeq \Map(M,\,N)$.
Applying this objectwise in the arrow category, any morphism from a Smith
ideal to a periodic Smith ideal factors, uniquely up to contractible choice,
through its periodization. By Proposition~\ref{periodic-app},
$\mathcal{P}^n(j)$ is periodic. Hence the canonical morphism
$j \to \mathcal{P}^n(j)$ factors through the periodization
$Pj: PI \to PA$. The nilpotent approximation functor $\mathcal{P}^n: \CAlg(\mathrm{Ar}^\Box(\mathbf{Spt}^\Sigma_{\mathbb{P}^1_+}(S))) \to \CAlg(\mathrm{Ar}^\Box(\mathbf{Spt}^\Sigma_{\mathbb{P}^1_+}(S)))$ of degree $n+1$ is homotopically idempotent. Therefore $\mathcal{P}^n(j)$ is a retract of $\mathcal{P}^n(Pj)$. 
\qed 
\end{proof}

\subsection{A motivic application: a \texorpdfstring{$K$}{K}-theoretic approximation of algebraic cobordism}

We conclude with a geometric application of the general theory developed above. Our guiding idea is that, rather than asking directly which $K$-theoretic properties should hold for algebraic cobordism itself, one can study the nilpotent approximation associated with the canonical morphism $\mathbf{MGL}\to\K$. This produces a natural $K$-theoretic approximation of algebraic cobordism. Its homotopical completeness will be a formal consequence of the general results proved above, while Bott periodicity and mod-$l$ rigidity will follow by combining the same formalism with standard geometric input from motivic homotopy theory. We recall the definitions of the homotopy $K$-theory spectrum and algebraic cobordism, following Voevodsky~\cite{VoeH}. For any smooth $S$-scheme $X$ and integers $n,\,r \ge 0$, let $\mathrm{Gr}_S(n,\,r)(X)$ denote the set of rank-$n$ subvector bundles of $\mathcal{O}_X^{n+r}$. Let $\mathbf{Gr}_S(n,\,r)$ denote the nerve of the Nisnevich sheafification of the presheaf $\mathrm{Gr}_S(n,\,r)(-)$. This defines a motivic space $\mathbf{Gr}_S(n,\,r)$. For $n \ge 0$, let $\mathrm{BGL}_{n\,/S}$ denote the motivic space $\varinjlim_r \mathbf{Gr}_S(n,\,r)$. We write $\proj^\infty = \mathrm{BGL}_{1\,/S}$, $\mathrm{BGL}_{/S}=\varinjlim_n \mathrm{BGL}_{n\,/S}$, and $\K=\Z \times \mathrm{BGL}_{/S}$, and call $\K$ the {\it (homotopy) $K$-theory spectrum}.

The algebraic cobordism spectrum $\mathbf{MGL}$ is the motivic spectrum determined by the Thom spaces of the universal vector bundles $V_n$ over $\mathbf{BGL}_n$ for all $n \ge 0$. It is known that the motivic spectrum $\mathbf{MGL}$ has a canonical $\E_\infty$-ring structure~\cite{Hu}. The motivic $\E_\infty$-ring $\mathbf{MGL}$ has the following universal property:
\begin{theorem}[Panin, Pimenov and R\"ondigs~\cite{PPR} Theorem 2.7]
\label{PPR}
Let $R$ be a motivic $\E_\infty$-ring. Then the set of monoidal maps $\mathbf{MGL} \to R$ is naturally isomorphic to the set of orientations on $R$. \qed
\end{theorem}

Let $\beta: \mathbb{P}^1 \to \mathbb{P}^\infty$ denote the Bott element.
Then it is known that $\beta \wedge -: \mathbb{P}^1_+ \wedge \K \to \K$ induces a weak equivalence $\K \to \Omega^1_{\proj^1}\K$. This weak equivalence is called {\it Bott periodicity}. Gepner--Snaith~\cite[Corollary 5.8]{GeSn} proved that $\K$ has a canonical $\E_\infty$-ring structure; therefore, the $K$-theory spectrum is a periodic motivic $\E_\infty$-ring.

We recall the following splitting principle:
\begin{theorem}[\cite{GeSn} Theorem 4.17]
\label{K-proj}
 Let $\beta: \mathbb{P}^1 \to \mathbb{P}^\infty$ denote the Bott element. Then the inclusion $\varphi:\mathbb{P}^\infty \to \mathrm{BGL}$ induces a motivic equivalence
\[
   \Sigma^\infty_+ \varphi_+: \Sigma_+^\infty \mathbb{P}_+^\infty [\beta^{-1}] \to \K, 
\] 
where $ \Sigma_+^\infty \mathbb{P}_+^\infty [\beta^{-1}]$ is the homotopy colimit of the diagram
\[
  \Sigma_+^\infty \mathbb{P}_+^\infty  \overset{\beta}{\to} \Omega_+^1 \Sigma_+^\infty \mathbb{P}_+^\infty\overset{\beta}{\to} \Omega_+^2 \Sigma_+^\infty \mathbb{P}_+^\infty  \overset{\beta}{\to} \cdots
\]
\qed
\end{theorem}

\begin{proposition}
\label{Bott-periodicity}
For each $n \ge 0$, the nilpotent approximation $\mathcal{P}^n_\K(\mathbf{MGL})$ of degree $n+1$ satisfies Bott periodicity; therefore so does $\Lambda_\K(\mathbf{MGL})$.
\end{proposition}
\begin{proof}
This follows from Proposition~\ref{split-P} and Bott periodicity in homotopy $K$-theory. \qed
\end{proof}

By Gepner--Snaith~\cite[Corollary 3.10]{GeSn}, the canonical map
\[
  \theta: P \mathbf{MGL} \to \Sigma_+^\infty \mathbf{BGL}[\beta^{-1}]
\]
obtained from the universal property in Theorem~\ref{PPR} is a weak equivalence of periodic motivic $\E_\infty$-algebras. In other words, if $P$ denotes periodization, then $P\mathbf{MGL}\simeq P\Sigma_+^\infty\mathbf{BGL}$. Moreover, by \cite[Proposition 5.10]{GeSn}, one has $\K\simeq P\Sigma_+^\infty B\mathbb{G}_m$. Therefore, the determinant map $\det:\mathbf{GL}\to\mathbb{G}_m$ induces a homotopically epimorphic morphism $P\mathbf{MGL}\to\K$.
\begin{proposition}
\label{finite-PMGL}
The motivic $\E_\infty$-ring $\K$ is a weakly compact $P\mathbf{MGL}$-module.
\end{proposition}

\begin{proof}
As a $P\mathbf{MGL}$-module, $P\mathbf{MGL}$ is the free module of rank one, hence compact. As explained above, the determinant map $\det:\mathbf{GL}\to\mathbb{G}_m$ induces a homotopically epimorphic morphism $P\mathbf{MGL}\to\K$. Applying Lemma~\ref{w-compact}, we conclude that $\K$ is weakly compact as a $P\mathbf{MGL}$-module. \qed
\end{proof}

The homotopical completeness of the motivic approximation is now a formal consequence of Proposition~\ref{finite-PMGL}, Proposition~\ref{compact-comp}, and the retract statement in Proposition~\ref{split-P}.
\begin{corollary}
\label{TheoremA}
The nilpotent approximations $\Lambda_\K(\mathbf{MGL})$ and $\Lambda_\K (P \mathbf{MGL})$ by algebraic $K$-theory are homotopically complete. 
\end{corollary}
\begin{proof}
Set $X=\Lambda_\K(\mathbf{MGL})$ and $Y=\Lambda_\K(P\mathbf{MGL})$. By Proposition~\ref{finite-PMGL} and Proposition~\ref{compact-comp}, $Y$ is homotopically complete. Moreover, Proposition~\ref{split-P} yields, for each $n \ge 0$, morphisms
\[
  \mathcal{P}^n_\K(\mathbf{MGL}) \longrightarrow \mathcal{P}^n_\K(P\mathbf{MGL}) \longrightarrow \mathcal{P}^n_\K(\mathbf{MGL})
\]
whose composite is the identity. Passing to the homotopy limit over $n$, we obtain morphisms $i:X \to Y$ and $r:Y \to X$ with $r \circ i=\mathrm{id}_X$. Let $\eta_X:X \to \Lambda(X)$ and $\eta_Y:Y \to \Lambda(Y)$ denote the canonical maps. By naturality,
\[
  \eta_X = \Lambda(r) \circ \eta_Y \circ i,
\]
so $\eta_X$ is a retract of $\eta_Y$. Since $Y$ is homotopically complete, $\eta_Y$ is a weak equivalence; hence so is $\eta_X$. Therefore $X$ is homotopically complete as well. \qed
\end{proof}

To relate this motivic approximation to the usual rigidity statement for homotopy $K$-theory, we first establish the following lifting theorem.
\begin{theorem}
\label{theorem-lift}
Let $f: X \to Y$ be a morphism between compact motivic spaces such that
\[
   f^*: \Map(\Sigma^\infty_+ Y,\,\K) \to \Map(\Sigma^\infty_+ X,\,\K)
\]
is a weak equivalence of motivic spectra. Then this weak equivalence lifts to the morphism
\[
  \mathcal{P}^n(f^*): \Map(Y,\, \mathcal{P}^n_\K(\mathbf{MGL})) \to \Map(X,\, \mathcal{P}^n_\K(\mathbf{MGL}))
\]
of nilpotent approximations of degree $n+1$ for each $n \ge 0$.
\end{theorem}
\begin{proof}
Let $j:I \to \mathbf{MGL}$ denote the homotopy kernel of the canonical
map $\mathbf{MGL} \to \K$. By Lemma~\ref{graded} (applied to $\mathcal{M}=\mathbf{Spt}^\Sigma_{\mathbb{P}^1_+}(S)$), for each $n\ge 1$ one has a homotopy fiber sequence
\[
  L_1(I^n/I^{n+1})  \to \mathcal{P}^{n}(j)  \to  \mathcal{P}^{n-1}(j)
\]
in the stable model category $\mathrm{Ar}^\Box(\mathbf{Spt}^\Sigma_{\mathbb{P}^1_+}(S))$.

By Lemma~\ref{n/n+1}, together with $\mathrm{Coker}(j)\simeq \K$, one has a canonical weak equivalence
\[
\K\otimes_{\mathbf{MGL}} I^n \to I^n/I^{n+1}.
\]
Hence $L_1(I^n/I^{n+1})$ has a canonical $\K$-module structure. Let $\mathcal{T}$ be the full subcategory of the homotopy category of $\K$-modules spanned by those $\K$-modules $M$ for which
\[
  f^*(M): \Map(Y,\,M) \to \Map(X,\,M)
\]
is a weak equivalence. Since $\Map(Y,-)$ and $\Map(X,-)$ are exact functors, $\mathcal{T}$ is thick. Because $X$ and $Y$ are compact motivic spaces, these functors commute with filtered homotopy colimits, so $\mathcal{T}$ is also closed under filtered homotopy colimits. The hypothesis gives $\K\in\mathcal{T}$.

The stable $\infty$-category underlying the model category of $\K$-modules is compactly generated by the free $\K$-module of rank one: by Lurie~\cite[Proposition~7.2.4.2]{HA}, the module category over an $E_1$-ring is compactly generated and its compact objects are precisely the perfect modules. Thus every $\K$-module is a filtered homotopy colimit of compact $\K$-modules, and the compact $\K$-modules form the thick subcategory generated by $\K$. Since $L_1(I^n/I^{n+1})$ is a $\K$-module, it is a filtered homotopy colimit of such compact modules; hence $L_1(I^n/I^{n+1})$ belongs to $\mathcal{T}$. Therefore
\[
  f^*(L_1(I^n/I^{n+1})): \Map(Y,\, L_1(I^n/I^{n+1})) \to \Map(X,\, L_1(I^n/I^{n+1}))
\]
is a weak equivalence for each $n \ge 0$.

For $n=0$, the assertion is exactly the hypothesis, equivalently the case $\mathcal{P}^0_\K(\mathbf{MGL})\simeq \K$. Assume that the assertion holds for $n-1$. Applying $\Map(Y,-)$ and $\Map(X,-)$ to the above homotopy fiber sequence, we obtain a morphism of homotopy fiber sequences whose left and right vertical maps are weak equivalences. Hence the middle vertical map is a weak equivalence. Therefore, by induction on $n$, for each $n \ge 0$, the induced morphism
\[
  \mathcal{P}^n(f^*): \Map(Y,\, \mathcal{P}^n_\K(\mathbf{MGL})) \to \Map(X,\, \mathcal{P}^n_\K(\mathbf{MGL}))
\]
is a weak equivalence. \qed
\end{proof}

\begin{corollary}
\label{theoremB}
Let $f: X \to Y$ be a morphism between compact motivic spaces. If
\[
   f^*: \Map(\Sigma^\infty_+ Y,\,\K) \to \Map(\Sigma^\infty_+ X,\,\K)
\]
is a weak equivalence of motivic spectra, then the induced maps
\[
  \Lambda_\K(f^*): \Map(Y,\, \Lambda_\K(\mathbf{MGL})) \to \Map(X,\, \Lambda_\K(\mathbf{MGL}))
\]
and
\[
  \mathcal{P}^* \Lambda_\K(f^*): \Map(Y,\, \Lambda_\K(P \mathbf{MGL})) \to \Map(X,\, \Lambda_\K(P\mathbf{MGL}))
\]
are weak equivalences.
\end{corollary}

\begin{proof}
By Theorem~\ref{theorem-lift}, the induced morphism of towers
\[
\{\Map(Y,\,\mathcal{P}^n_\K(\mathbf{MGL}))\}_{n\ge 0}\to \{\Map(X,\,\mathcal{P}^n_\K(\mathbf{MGL}))\}_{n\ge 0}
\]
is levelwise a weak equivalence. Since homotopy limits preserve levelwise weak equivalences, passing to homotopy limits yields
\[
\Lambda_\K(f^*): \Map(Y,\, \Lambda_\K(\mathbf{MGL})) \to \Map(X,\, \Lambda_\K(\mathbf{MGL})).
\]

Let $j_P:J\to P\mathbf{MGL}$ denote the homotopy kernel of the canonical map $P\mathbf{MGL}\to\K$. Repeating the proof of Theorem~\ref{theorem-lift} with $j$ replaced by $j_P$, we obtain a levelwise weak equivalence of towers
\[
\{\Map(Y,\,\mathcal{P}^n_\K(P\mathbf{MGL}))\}_{n\ge 0}\to \{\Map(X,\,\mathcal{P}^n_\K(P\mathbf{MGL}))\}_{n\ge 0}
\]
Applying the same observation to homotopy limits, we obtain
\[
  \mathcal{P}^* \Lambda_\K(f^*): \Map(Y,\, \Lambda_\K(P \mathbf{MGL})) \to \Map(X,\, \Lambda_\K(P\mathbf{MGL}))
\]
as a weak equivalence. \qed
\end{proof}

\begin{theorem}
\label{rigidity}
Let $A$ be a commutative ring and $I$ an ideal such that $(A,\,I)$ is a Henselian pair and the residue ring $A/I$ has positive characteristic $p$. Let $l$ be an integer invertible in $A/I$. Then the closed immersion
\[
  i:\Spec{A/I} \to \Spec{A}
\]
induces a weak equivalence
\[
 i^*:  \Lambda_{\K/l}(\mathbf{MGL}/l)(\Spec{A}) \to \Lambda_{\K/l}(\mathbf{MGL}/l)(\Spec{A/I}), 
\]
where $\mathbf{MGL}/l$ and $\K/l$ denote the mod-$l$ algebraic cobordism and $K$-theory, respectively.
\end{theorem}
\begin{proof}
Let $i:\Spec(A/I)\to\Spec(A)$ be the closed immersion. Since $(A,\,I)$ is a Henselian pair and $l$ is invertible in $A/I$, Gabber rigidity for homotopy $K$-theory~\cite{Gabber-rigid} yields a weak equivalence
\[
  i^*: \Map(\Sigma^\infty_+\Spec(A),\,\K/l) \xrightarrow{\sim} \Map(\Sigma^\infty_+\Spec(A/I),\,\K/l).
\]
Replacing $\mathbf{MGL}\to\K$ by the quotient morphism $\mathbf{MGL}/l\to\K/l$ in the proofs of Theorem~\ref{theorem-lift} and Corollary~\ref{theoremB}, we obtain the corresponding mod-$l$ statement. Applying it to $f=i$, we obtain
\[
 i^*:  \Lambda_{\K/l}(\mathbf{MGL}/l)(\Spec{A}) \to \Lambda_{\K/l}(\mathbf{MGL}/l)(\Spec{A/I}).
\]
This is the desired weak equivalence. \qed
\end{proof}

\subsubsection*{Acknowledgements}  
The author used Google’s Gemini-pro 3.1 as a conversational research aid for brainstorming and refining mathematical formulations, and OpenAI’s Prism (GPT-5.2) for editorial and expository assistance, including structural refinement. These tools were used to improve the manuscript’s clarity and to indicate where abbreviated arguments warranted fuller exposition. All results, proofs, and final formulations were reviewed, verified, and approved by the author, who bears sole responsibility for the content of this work.
\begingroup
\sloppy
\bibliographystyle{alphadin}
\bibliography{bibkato}
\endgroup
\nocite{kato2023mathematics}
\nocite{zbMATH03811819}
\nocite{zbMATH07335469}
\nocite{zbMATH05706074}
\nocite{zbMATH02095715}
\nocite{zbMATH01787461}
\nocite{Hirschhorn}
\nocite{Hoveybook}
\nocite{HT}
\end{document}